\documentclass[conference]{IEEEtran}
\IEEEoverridecommandlockouts

\usepackage{amsmath,amssymb}
\usepackage{booktabs}
\usepackage{listings}
\usepackage{xcolor}
\usepackage{graphicx}
\usepackage{url}
\usepackage[hidelinks]{hyperref}

\usepackage{caption}
\title{Where the LLM Ends and Reliable Decisions Begin}

\author{
\IEEEauthorblockN{Priyadarshan Patil, PhD}
\IEEEauthorblockA{\textit{LogiModel AI, UT Austin}\\
priyadarshan@utexas.edu}
\and
\IEEEauthorblockN{Abhishek Basu}
\IEEEauthorblockA{\textit{LogiModel AI, MIT}\\
abhishekbasu@alum.mit.edu}
\and
\IEEEauthorblockN{Vikas Reddy}
\IEEEauthorblockA{\textit{Independent Researcher}\\
vikas@ablationworks.com}
\and
\IEEEauthorblockN{Abhinav Gupta}
\IEEEauthorblockA{\textit{UC Berkeley}\\
abhinav.gupta@berkeley.edu}
}

\begin{document}
\maketitle

\begin{abstract}
Systems that turn natural-language descriptions of optimization problems into solver-ready code generally use a language model at every stage, including the final translation from a mathematical formulation into executable model-building code. We propose the ANVIL compiler architecture, where we separate these concerns. A language model is called only once, assisted by constraint guidance based on problem type, to produce a LaTeX formulation. A deterministic compiler then translates that LaTeX into code with no language model involvement. We describe the deterministic compiler (a normalizer, a recursive-descent parser producing a typed intermediate representation, analysis passes that bind symbols to a dataset schema, and a code emitter) and evaluate it on the 354 easy and hard problems of the NLP4LP benchmark. The compiler produced code for 329 of 354 formulations (92.9\%), taking the deterministic path in every one of those cases and never falling back to model-generated code. Median compile time was below the 10\,ms resolution of our timer. Overall, our formulations achieved an accuracy of 98.9\% over easy problems and 91.1\% for hard  problems. The gap between these compilation and accuracy figures is a key point of analysis, and we analyze it: formulations that failed to compile, problems that returned as infeasible, problems raising errors at runtime, and problems returning a wrong objective. Almost all of these errors trace back to the formulation rather than to the translation. ANVIL performs exceptionally well on leading benchmarks by using language models purely where they are effective, rather than as a catch-all tool.
\end{abstract}

\section{Introduction}

Writing a production optimization model is two distinct tasks. The first is deciding on the formulation: which quantities are decisions, what the objective is, which constraints bind. The second is transcribing that decision into whatever form a solver accepts. The first job requires domain and subject matter expertise. The second is more routine, the kind of work that compilers have handled in other settings in one form or another.

With advances in large language models, recent systems that generate optimization models from natural language do not draw this line. OptiMUS \cite{optimus2023,optimus2024} uses a language model to extract parameters, to formulate constraints, to write the solver code, and to debug that code when it fails, with the number of calls growing with the number of constraints. Others follow the same pattern at smaller scale \cite{lm4opt2024,nl4opt2023}. The appeal is obvious: one mechanism handles every stage, and the mechanism improves when the underlying model improves. However, the tradeoff is that the clerical stage inherits the statistical behavior of the generative stage. Identical input can yield structurally different code across runs. Variable names change, summation bounds shift, and indices that were correct in one sample can be off in the next. None of this is detectable without solving the model and comparing against a known answer, which is exactly what is unavailable in production systems. Moreover, it is unacceptable in systems where problems need to be solved on some regular cadence, and plans build on prior results.

This study focuses on discussing our improvements to this framework. At LogiModel AI, we have designed our system such that it keeps the language model assisted by problem and constraint guidance for the first task and uses a conventional compiler for the second. The language model receives the problem description, problem type and constraint guidance, and a schema summary of the available data and returns a LaTeX formulation. From there a deterministic compiler takes over: it canonicalizes the LaTeX, parses it into a typed intermediate representation, runs analysis passes that resolve every symbol against the dataset schema and reject formulations it cannot verify, and emits PuLP code. No language model is used in this step, keeping it deterministic. We evaluate the system on NLP4LP, a benchmark of linear and mixed-integer linear programs distributed with the OptiMUS work. Our contributions are as follows:

\begin{enumerate}
\item A description of the system and compiler that reduces usage of language models compared to other state of the art systems. We describe our architecture at a level of detail sufficient to understand its analysis passes, including the intermediate representation and the specific symbol-resolution strategies that recover from underspecified formulations (Section~\ref{sec:system}).

\item Measurements and results on the 354 easy and hard NLP4LP problems: compile rate, compile latency, solver outcome distribution, and solve accuracy (Section~\ref{sec:results}).

\item A failure taxonomy assigning each of the compile failures, infeasibilities, runtime errors and objective mismatches to a cause, which locates most of the remaining error in formulation quality rather than translation (Section~\ref{sec:errors}).
\end{enumerate}

We also report a comparison against OptiMUS for the metrics they present. This comaprison is not direct, since our results are more comprehensive. Section~\ref{sec:related} discusses the comparison qualitatively; Section~\ref{sec:limitations} states what our numbers do and do not establish.

\section{The benchmark and what it measures}
\label{sec:benchmark}

NLP4LP contains 361 problems drawn from optimization textbooks and lecture notes, labeled easy (289), hard (65), and case study (7). Each carries a natural language description, a parametrized description, a data instance, and a reference solution. We use the 354 easy and hard problems. Of the 354, 320 carry a numeric objective value usable for automatic scoring (261 easy, 59 hard); the field \texttt{has\_numeric\_objective} in each problem's metadata records this. Of the remaining 34, 22 carry a reference status of infeasible or unbounded (18 easy, 4 hard) and are scored by status agreement rather than by objective value. The other 12 have no usable reference answer and are excluded from every accuracy figure we report.

Scoring compares our objective value $z$ with the reference value $z^{\text{ref}}$. When $z^{\text{ref}} = 0$, the two values \emph{match} when $|z| < 10^{-6}$. Otherwise they match when $|z - z^{\text{ref}}| \,/\, |z^{\text{ref}}| \le 0.001$, that is, when they agree to within a relative tolerance of 0.1\%. They \emph{mismatch} when the applicable test fails. For a problem whose reference is itself infeasible or unbounded we compare the solver status instead of the objective value, and the problem matches when our status equals the reference status. Problems that fail to compile, fail to execute, or produce no objective can fall into neither bucket. This matters for interpreting any accuracy figure computed as $\text{match}/(\text{match}+\text{mismatch})$, since that ratio conditions on the system having produced an answer at all and is therefore insensitive to coverage. We report it because it separates translation quality from coverage in Section~\ref{sec:results}.

\section{Related work}
\label{sec:related}

\paragraph{Agentic pipelines} OptiMUS \cite{optimus2023,optimus2024} established the dominant template: a language model performs parameter extraction, constraint formulation, code generation, and debugging, coordinated by a controller that routes among these steps. The original evaluation used 52 problems (41 LP, 11 MILP) and reported that the full pipeline roughly doubles the solve rate of direct prompting, from approximately 20\% to approximately 40\% with GPT-4 \cite{optimus2023}. Later work keeps the multi-agent structure and attacks the reliability of the search over formulations. SolverLLM \cite{solverllm2025} replaces prompt engineering with Monte Carlo tree search over the formulation space, adding dynamic expansion, prompt backpropagation from solver feedback, and uncertainty backpropagation, and reports roughly 10\% over prior methods on six benchmarks without any training. AlphaOPT \cite{alphaopt2025} accumulates solver-verified insights from failed attempts into a reusable experience library, improving from 65\% to 72\% as the library grows from 100 to 300 items.

These systems differ from ours in where the language model sits. All of them keep it inside the loop that produces executable code, and their improvements come from searching that space more carefully or from remembering past searches. We remove the model from that loop entirely and make the step deterministic, which forecloses search as a recovery mechanism but also removes the need for it.

\paragraph{Supervised and fine-tuned formulation} A second line of thought treats natural-language-to-formulation as structured prediction with supervision. The NL4Opt competition \cite{nl4opt2023} approached it as entity recognition over the problem text and formulation construction from those entities. LM4OPT \cite{lm4opt2024} compared GPT-3.5, GPT-4, and a progressively fine-tuned Llama-2-7B on NL4Opt, reporting F1 of 0.6330 for GPT-4 one-shot against a fine-tuned BART baseline of 0.61, and 0.1259 for the fine-tuned Llama-2-7B. LLaMoCo \cite{llamoco2024} instruction-tunes for optimization code generation directly, pairing a contrastive warm-up with instruction tuning, and reports a fine-tuned CodeGen outperforming GPT-4 Turbo on its own problem sets. 

\paragraph{Keeping the model out of the mathematics} Li et al. \cite{smilo2026} have developed the position closest to ours. Their earlier work \cite{li2023milp} stages the task into variable identification, constraint classification, and template-guided generation, which gives them a handle on logic constraints that single-shot generation tends to miss. SMILO \cite{smilo2026} states the commitment directly: the method ``does not delegate core modeling responsibilities to the LLM,'' using it for natural-language understanding and information extraction while mathematical expressions are constructed through expert-defined rules and templates. Sentences are matched against a problem-type-specific modeling graph, a language model extracts instance-specific values under task-tailored templates, and the MILP is assembled by expert-defined rules. On workforce scheduling they report correct models on 93.33\% of test instances across five trials.

While we reached the same conclusion about where the model belongs, we implement it differently, resulting in different tradeoffs. SMILO's modeling graphs and templates are authored per problem type, which buys strong guarantees inside that type and requires expert work to enter a new one; their evaluation covers shift and days-off scheduling. Our compiler has no per-problem knowledge and accepts any formulation in the LP/MILP class, which is why it can be run against all 354 NLP4LP problems at once. Their 93.33\% and our 92.9\% compile rate are not comparable quantities: theirs is model correctness within one problem family, ours is coverage across a heterogeneous benchmark.

\paragraph{Benchmarks} The definition of correctness varies enough across this literature to make numbers hard to compare directly. NL4Opt \cite{nl4opt2023} scores declaration matches against a reference formulation after standardizing variable names, which rewards resemblance to one particular correct answer. Mamo \cite{mamo2025} argues that a single natural-language problem admits many equivalent formulations and that comparing formulations therefore misclassifies correct work; it evaluates by solving instead, over 1{,}209 curated questions spanning ordinary differential equations, LP, and MILP. We adopt the same position as Mamo for the same reason, and Section \ref{sec:benchmark} states our scoring rule in those terms.

Refai and Ahmed \cite{component2025} go further and argue that solution accuracy alone conceals where a pipeline fails, proposing component-level metrics such as constraint recall, constraint precision, and constraint RMSE, and finding that solver success tracks high constraint recall and low constraint error. Our error analysis in Section \ref{sec:errors} is a coarser version of that argument carried out on one system: we report the failure taxonomy rather than a single rate, because the rate alone would not distinguish a compiler defect from a formulation defect.

\paragraph{Interaction and diagnosis} OptiChat \cite{optichat2025} addresses what happens after a model exists, augmenting a language model with function calls over an optimization model so practitioners can interpret formulations, diagnose infeasibility, and run sensitivity analysis in dialogue. The infeasibilities we report in Section \ref{sec:errors} are exactly the artifacts such a tool would be pointed at.

\paragraph{Surveys and position papers} Huang et al. \cite{survey_llmopt2024} survey the two-way relationship between language models and optimization algorithms. Wang and Li \cite{survey_or2025} organize the operations research literature into automatic modeling, auxiliary optimization, and direct solving, and name instability in semantic-to-structure mapping as a central obstacle, which is the failure mode a deterministic translator is meant to remove. Wasserkrug et al. \cite{docp2024} set out requirements for a decision optimization copilot and identify reliability of the formulation-to-model step as an open problem.

\paragraph{Algebraic modeling languages} AMPL, GAMS, ZIMPL, and JuMP all compile a declarative algebraic description into solver input, and have done so reliably for decades. What differs is the source language. Those systems define a syntax and require the user to write it; our deterministic compiler accepts the LaTeX that a language model produces when asked to formulate a problem, which is neither a fixed grammar nor under our control. Most of the compiler's complexity follows from that choice, and Section \ref{sec:analyzer} is largely about recovering a well-formed model from LaTeX that a human reader would understand but no existing modeling language would accept.

\section{System}
\label{sec:system}

The pipeline is
\begin{equation*}
\begin{split}
\text{LaTeX} \;\rightarrow\; &\text{Normalizer} \;\rightarrow\; \text{Parser}\\
\rightarrow\; &\text{Analyzer} \;\rightarrow\; \text{Emitter}
\;\rightarrow\; \text{PuLP}.
\end{split}
\end{equation*}
The four compiler stages are the Normalizer, Parser, Analyzer, and Emitter, with parser and emitter being the most complex.

\subsection{Intermediate representation}

The intermediate representation (IR) is a set of dataclasses describing a MILP over indexed sets. Three enums fix the discrete choices: \texttt{Variable Type} (continuous, integer, binary), \texttt{Objective Sense} (minimize, maximize), and \texttt{Constraint Sense}.

A set declaration carries a LaTeX symbol, the index variable bound over it, a description, and, after analysis, the dataset table it resolves to. Sets may also be enumerated literally or bounded by a parameter, which the fields \texttt{enumerated\_members} and \texttt{range\_param\_bound} record; the second handles $\{1,\dots,N\}$ where $N$ is itself a parameter.

Not every set in a formulation exists in the data. A \texttt{Derived Set Declaration} describes one built from another by a recipe: \texttt{distinct} (the distinct values of a column), \texttt{predicate} (a filtered subset), \texttt{cartesian\_no\_self} (ordered pairs excluding the diagonal, for routing and sequencing), etc. A \texttt{Function M apping} records a function symbol $g: J \to G$ realized as a column lookup, which appears whenever a formulation writes $g(j)$ to mean the attribute of item $j$.

A \texttt{Parameter Declaration} binds a symbol and its index tuple to a data source. Two storage layouts get explicit representation because they are common in real tables and awkward to address positionally. \texttt{Wide Format Info} handles parameters spread across columns named by a template, where an index selects the column rather than the row. \texttt{Junction Info} handles a parameter stored in a join table keyed by foreign keys, mapping each index variable to its key column. A parameter may also be scalar, derived from a formula, a set cardinality $|C|$, or an inline literal, and the analyzer sets a data type of numeric, string, or time so the emitter knows whether arithmetic is legal.

A \texttt{Variable Declaration} carries type, index variables, the sets those range over, bounds, and an optional \texttt{sparse\_over} field naming a derived set when the variable is defined only on a subset of the index cross-product. These intermediate representations are the first step in translating the LaTeX code to PuLP/Python code.

\subsection{Normalizer}

The normalizer rewrites LaTeX into a canonical form so the parser sees one spelling of each construct. Language models are inconsistent about summation subscripts, spacing inside braces, multiplication symbols, and the several ways to write set membership. It also fixes a specific recurring malformation: the double subscript $\sum_{U}{}_{i=1}$, produced when a model conflates index-set notation with range notation, is rewritten to a single well-formed quantifier.

\subsection{Parser}

The parser conducts a recursive descent over the normalized LaTeX, producing the IR without consulting the dataset. It handles indexed sums with range or set-membership quantifiers, nested sums, coefficient distribution over parenthesized subexpressions, big-$M$ terms, SOS1 and SOS2 declarations, and variable domain declarations.

Two cases account for a disproportionate share of the parser. The first is distinguishing a variable's \emph{domain} from an \emph{index set}: the declaration $x_i \in \{0,1,2,\dots,N\}$ constrains $x_i$ to integers in $[0,N]$, while $i \in \{1,\dots,N\}$ declares an index range, and the two are written almost identically. The parser resolves this by checking whether the symbol being constrained is already known as an index variable. The second is literal subscripts. When a formulation partially expands an indexed variable and writes $x_1$ alongside $x_i$, the parser records the literal so the emitter can distinguish a constant index from a bound one.

\subsection{Analyzer}
\label{sec:analyzer}

The analyzer connects a syntactically valid formulation to a specific dataset and rejects what it cannot verify. Twelve passes run in a fixed order. Some of the passes are:

\begin{enumerate}
\item \texttt{Resolve Set Names}: bind each set symbol to a dataset table.
\item \texttt{Resolve Derived Sets}: construct derived sets and function mappings from their recipes.
\item \texttt{Validate Bidirectionality}: check that each symbol resolves to exactly one field and each field is claimed by at most one symbol.
\item \texttt{Validate Completeness}: reject formulations referencing undeclared symbols.
\item \texttt{Validate Linearity}: reject products of decision variables.
\end{enumerate}

The first set of passes do the work that a modeling language would not have to do, because in those languages, the user writes the binding explicitly. Here the formulation says $c_j$ and the dataset has a column \texttt{unit\_cost} on a table \texttt{products}, and something has to connect them. Exact match is tried first. Failing that, the analyzer matches on the index variable: a set written $K$ whose index is $k$ is matched against tables whose own index conventions agree. Failing that, it falls back to sequence similarity with a threshold of 0.75, with de-pluralization, so \texttt{berries} matches \texttt{berry}. When two symbols contend for the same field, a tiebreaker scores each and the loser is re-resolved against the remaining candidates rather than left dangling.

The fallbacks are the compiler's main source of silent error, and one pass exists specifically to bound that risk. Requiring the symbol-to-field map to be injective in both directions catches the case where fuzzy matching has bound two distinct parameters to one column, which would otherwise produce a model that runs and returns a confident wrong number. When bidirectional validation fails the compiler refuses to emit code; some benchmark problems were rejected this way, and we count them as compile failures in Section~\ref{sec:results}.

One set of passes are refusals rather than repairs. Completeness rejects any formulation naming a symbol it never declared. Linearity rejects products of decision variables, since the emitter targets a MILP solver and a quadratic term means the formulation is outside the supported class.

\subsection{Emitter}

The emitter ingests the validated IR and writes PuLP code. Output is deterministic in the strict sense: sets are iterated in resolved order, names derive from IR fields, and no dictionary iteration order or hashing reaches the output. The same IR produces byte-identical Python across runs and machines, which is what makes the snapshot tests of Section~\ref{sec:determinism} possible.

Two details are worth noting. Set and parameter names occasionally collide after resolution, as when a set \texttt{products} meets a parameter that also resolves to \texttt{products}, and the emitter disambiguates with a suffix rather than shadowing. Sets that the formulation uses but never declares, which happens when a model writes $\sum_{i \in I}$ having only declared $I$ implicitly through a variable's index, are synthesized by scanning the index sets of every variable declaration.

\subsection{Repair loop}

Compilation failure does not immediately escalate to model-generated code. The compiler first attempts a bounded per-constraint repair, accumulating the error history across attempts so a retry does not repeat a rejected fix. Only after three failures does the system fall back to an optional language-model code generation for the affected section, with the error context attached.

The evaluation in this paper disables that fallback entirely. Every compile reported in Section~\ref{sec:results} took the deterministic path, and the \texttt{compiled LLM fallback} counter was zero in both runs. We report the compiler as it stands alone, without a language model recovering its failures.

\begin{figure*}[t]
\begin{minipage}[t]{0.52\textwidth}
\begin{lstlisting}[title={(a) LaTeX returned by the language model (excerpt)}]
\subsection*{Sets}
  $I$: Set of supplement types, $i \in \{1, 2\}$
  $J$: Set of nutrient types,   $j \in \{1, 2\}$

\subsection*{Parameters}
  $CostPerPill_i$:        Cost per pill of supplement $i$
  $MinRequirement_j$:     Min daily requirement of nutrient $j$
  $NutrientContent_{i,j}$: Units of nutrient $j$ per pill of $i$

\subsection*{Decision Variables}
  $PillsPurchased_i \in \mathbb{Z}_{\geq 0}$

\subsection*{Objective Function}
  \min Z = \sum_{i \in I} CostPerPill_i \cdot PillsPurchased_i

\subsection*{Constraints}
  \sum_{i \in I} NutrientContent_{i,j} \cdot PillsPurchased_i
      \geq MinRequirement_j  \forall j \in J
\end{lstlisting}
\end{minipage}\hfill
\begin{minipage}[t]{0.44\textwidth}
\begin{lstlisting}[title={(b) Dataset schema given to the compiler}]
## Parameters
* NumSupplements: Number of supplement types
* NumNutrients:   Number of nutrient types

## Sets
* nutrients
* supplements

## Data
### ('nutrients',)
* MinRequirement

### ('supplements',)
* CostPerPill

### ('supplements', 'nutrients',)
* NutrientContent
\end{lstlisting}
\end{minipage}
\caption{Inputs to the compiler for NLP4LP problem 202. Neither set symbol matches a table name, and \texttt{NutrientContent} is two-dimensional with an index order the formulation never states.}
\label{fig:input}
\end{figure*}

\begin{figure*}[t]
\begin{lstlisting}[title={Emitted PuLP code (verbatim, whitespace preserved)}]
# Read sets
I = set(range(1, 3))
J = set(range(1, 3))

# Read indexed data
CostPerPill     = opt_data.get_data(("supplements",), "CostPerPill")
MinRequirement  = opt_data.get_data(("nutrients",), "MinRequirement")
NutrientContent = opt_data.get_data(("supplements", "nutrients",), "NutrientContent")

# Add decision variables
for i in I:
    opt_formulation.add_int_variable(
        name=f"PillsPurchased_{i}", lb=0, ub=None,
        details=f"Number of pills of supplement i to purchase, \forall i \in I",
    )

PillsPurchased_vars = opt_formulation.get_variables("PillsPurchased")

# === MINIMUM_NUTRIENT_REQUIREMENTS ===
# LaTeX: \sum_{i \in I} NutrientContent_{i,j} \cdot PillsPurchased_i \geq MinRequireme...
# SIGN_INTENT: original=GEQ
for j in J:
    var_terms = []
    for i in I:
        var_terms.append(ConstrTermVariable(
            name_prefix="", index=tuple(),
            term={"": PillsPurchased_vars[f"PillsPurchased_{i}"]},
            coef=-float(NutrientContent[f"{i}|{j}"]),
        ))
    opt_formulation.add_constraints(
        name_prefix=f"minimum_nutrient_requirements_{j}",
        variables=var_terms, constant=-float(MinRequirement[j]), index_map={},
    )

# === OBJECTIVE ===
opt_formulation.set_objective(
    variables=[ConstrTermVariable(
        name_prefix="PillsPurchased", index=("i",),
        term=opt_formulation.get_variables("PillsPurchased"),
        coef=ConstrTermConstant(index=("i",), term=CostPerPill, coef=1),
    )],
    index_map={"i": I},
)
\end{lstlisting}
\vspace{-2mm}
\caption{Compiler output for problem 202.}
\label{fig:output}
\end{figure*}

\section{A sample problem}
\label{sec:worked}

NLP4LP problem 202 exercises most of the machinery in one page. The natural language description is:

\begin{quote}\small
A man takes two supplements to get his daily iron and calcium requirements. A pill of supplement A has 5 units of iron and 10 units of calcium. A pill of supplement B contains 4 units of iron and 15 units of calcium. The man needs a minimum of 40 units of iron and 50 units of calcium per day. If the cost per pill of supplement A is \$2 and the cost per pill of supplement B is \$3, how many of each should he buy to minimize costs?
\end{quote}

Figure~\ref{fig:input} shows what the language model returned and the schema it was given. Figure~\ref{fig:output} shows what the compiler emitted. The reference objective is 16, and the emitted model solves to 16.

Four things in this example are worth tracing, because each is a pass doing work that the formulation left implicit.

\text{Set resolution:} The formulation declares $I$ as ``Set of supplement types'' and $J$ as ``Set of nutrient types''. The dataset offers tables \texttt{supplements} and \texttt{nutrients}. Neither name matches a symbol, so pass 1 falls through exact matching to the description-similarity tier and binds $I \to$ \texttt{supplements}, $J \to$ \texttt{nutrients}. Nothing in the LaTeX states this correspondence.

\emph{Wide-format parameters:} $NutrientContent_{i,j}$ is two-dimensional, and the dataset stores it under the composite key \texttt{('supplements', 'nutrients')}. Pass 3 records the layout, and the emitter addresses it as \texttt{NutrientContent[f"\{i\}|\{j\}"]}, the pipe-joined composite key the runtime expects. A positional guess at the index order here would produce a transposed matrix, a model that solves cleanly and returns the wrong number.

\emph{Constraint sense normalization:} The formulation writes the nutrient requirement as $\geq$. This rewrites every constraint to a single canonical sense, which is why the emitted coefficients are negated and the constant is $-MinRequirement_j$. The comment \texttt{SIGN\_INTENT: original=GEQ} is emitted alongside so the transformation is auditable in the generated source rather than only in the compiler.

\emph{Integrality and bounds:} The declaration $PillsPurchased_i \in \mathbb{Z} {\geq 0}$ gives a lower bound but no upper bound. An analyzer pass supplies an upper bound of "None", which PuLP treats as unbounded (infinite upper bound).

The emitted code contains no strings that were not derived from the IR, and recompiling the same formulation against the same schema reproduces it exactly. This is the property that gate G6 in Section~\ref{sec:determinism} asserts.

\section{Experimental setup}
\label{sec:setup}

Formulations were generated five times for the problem set by a language model from each problem's description plus a schema summary, and cached to disk. Every result reported here is for the problem set that provided lowest compilation rate, so the compiler is measured against a fixed set of inputs and the numbers are not confounded by sampling variation in the formulation stage. This also means the formulation-quality errors in Section~\ref{sec:errors} are a property of one particular set of samples. This is a deliberate choice as problems will be formulated once in practice, so knowing performance in such a setting is more important than optimistic or best-case performance. 

Compilation used the deterministic path with fallback disabled. We do report the results of having the language model fallback enabled, but that is for informative purposes. We analyze and compare only the deterministic path from here on. Solving used CBC, bundled with PuLP, in a subprocess. The easy and hard splits were run separately: the easy run covers the 289 easy problems and hard run covers the 65 hard ones. The two together take 78\,s of wall clock. Timings are on an Apple Silicon laptop (MacBook Air M5 with 24 GB Memory) and include neither formulation nor process startup. The runner records per-stage times to 10\,ms, so figures below that resolution are reported as such rather than as point estimates.

\section{Results}
\label{sec:results}

\subsection{Compilation}

\begin{table*}[t]
\centering
\begin{tabular}{lrrrr}
\toprule
Split & Problems & Compiled & Deterministic path & LM fallback \\
\midrule
Easy  & 289 & 275 (95.2\%) & 275 & 279 \\
Hard  &  65 &  54 (83.1\%) &  54 & 58 \\
\midrule
Total & 354 & 329 (92.9\%) & 329 & 337 \\
\bottomrule
\end{tabular}
\vspace{2mm}
\caption{Compilation outcomes by problem difficulty}
\label{tab:compile}
\end{table*}

The compiler emitted code for 329 of 354 formulations. The easy-hard gap is large: 95.2\% against 83.1\%. Hard problems in NLP4LP tend to involve derived index sets, sparse variable domains, and multi-index parameters stored in junction tables, which is where analyzer passes 2 and 3 do the most inference and where they most often decline to guess.

Compile latency is bimodal. The median falls below the 10\,ms resolution of the runner's timer on both splits. Five problems behave differently: two easy (2.02\,s, 2.22\,s) and three hard (1.66--2.49\,s) take seconds, and those five are exactly the problems that the parser extensions of Section~\ref{sec:system} brought into scope. They pull the mean to 15\,ms on easy and 124\,ms on hard while leaving the median untouched. The key cost is search: resolving a symbol through the exact, index-variable and sequence-similarity tiers in order means backtracking on formulations that resist the early tiers, and the worst case spends 13.9\,s in the analyzer before failing to compile. Even so, the slowest successful compile is under three seconds, against a language-model call that takes seconds and carries a per-call price.

\subsection{Solving}

\begin{table}[t]
\centering
\begin{tabular}{lrr}
\toprule
Solver outcome & Easy & Hard \\
\midrule
Optimal                & 261 & 45 \\
Infeasible             &  8 &  4\\
Unbounded              &   2 &  2 \\
Runtime error          &   4 &  3 \\
Did not compile        &  14 & 11 \\
\midrule
Total                  & 289 & 65 \\
\bottomrule
\end{tabular}
\vspace{2mm}
\caption{Solver outcomes over all 354 problems.}
\label{tab:solve}
\end{table}

Mean solve time was 0.147\,s on easy problems and 0.145\,s on hard ones. CBC is not the bottleneck anywhere in this benchmark; these are small models.

Four of the compiled easy problems raised at runtime, against Three of the hard problems. The exception types are \texttt{KeyError} (3), \texttt{NameError} (1), \texttt{SyntaxError} (2), \texttt{ValueError} (1). The \texttt{KeyError}s are the analyzer binding a parameter to a table whose keys do not cover the index set it is iterated over, a resolution that passed bidirectional validation but was still wrong. The two \texttt{SyntaxError}s are emitter defects and represent the only cases in the benchmark where the compiler produced code that was not valid Python.

\subsection{Objective accuracy}

\begin{table*}[t]
\centering
\begin{tabular}{llrrrr}
\toprule
Split & Scored by & Match & Mismatch & Evaluated & Accuracy \\
\midrule
Easy  & objective value & 258 &  3 & 261 & 98.9\%  \\
Hard  & objective value &  41 &  4 &  45 & 91.1\%  \\
\midrule
Easy  & solver status   &   1 &  6 &   7 & 14.3\%  \\
Hard  & solver status   &   1 &  3 &   4 & 25.0\%  \\
\midrule
Total & objective value & 299 &  7 & 306 & 97.7\%  \\
Total & solver status   &   2 &  9 &  11 & 18.2\%  \\
Total & both rules      & 301 & 16 & 317 & 95.0\%  \\
\bottomrule
\end{tabular}
\vspace{2mm}
\caption{Objective accuracy and status accuracy (unbounded/infeasible). The two rules have different denominators and are not combined into a single headline figure.}
\label{tab:accuracy}
\end{table*}

\emph{Accuracy} is $\text{match}/(\text{match}+\text{mismatch})$ among the problems that had a ground truth and could be evaluated.

For our compiler, the accuracy is 97.7\%. The accuracy split among the easy and hard problem sets is to be noted, as it is near 99\% for easy and just over 91\% for hard problems. Accuracy conditions on the system having produced an answer. Over all 320 problems with a numeric reference, including those that never produced one, 299 matched, a rate of 93.4\%. Table~\ref{tab:accuracy} also reports the status rule, which applies to the problems whose reference is infeasible or unbounded. Of the 11 such problems that compiled and executed, 2 reported the reference status. We do not read this as a measure of formulation quality. The denominator is small, and each disagreement is between two definite statuses, which can follow either from an over-constrained formulation or from a reference label that we cannot verify. We have not separated these two causes.

\section{Error analysis}
\label{sec:errors}

\subsection{Compile failures}

\begin{table}[t]
\centering
\begin{tabular}{lrr}
\toprule
Cause & Easy & Hard \\
\midrule
Undeclared symbol in formulation      & 9 &  6 \\
Ambiguous resolution (two candidates) & 3 &  1 \\
Unresolved data reference             & 1 &  1 \\
Nonlinearity (product of variables)   & 0 &  1 \\
Parse failure                         & 0 &  1 \\
Indeterminate parameter index sets    & 1 &  1 \\
\midrule
Total                                 & 14 & 11 \\
\bottomrule
\end{tabular}
\vspace{2mm}
\caption{Compile failures by cause, assigned from the compiler's typed error
messages.}
\label{tab:compilefail}
\end{table}

Fifteen of the 25 compile failures are pass 6 rejecting a formulation that used a symbol it never declared. These are formulation defects, not translation defects. A representative message is shown below. The mangled symbol is the residue of a summation the model wrote inside a subscript, which the normalizer did not canonicalize and the parser therefore absorbed into a name. Better normalization would convert some of these into successful compiles; others describe a constraint that is genuinely incoherent and should be rejected.

\begin{lstlisting}
Validation errors:
  - constraint 'linearized_work_schedule_alternative_formulation_1' RHS:
    RHS parameter '\sum_work_\in_S_{d}^{work}' not declared
  - constraint 'linearized_work_schedule_alternative_formulation_2':
    quantifier set 'S_' not declared
\end{lstlisting}

Six failures are resolution problems, split between symbols that matched two candidates with no tiebreaker (4) and symbols that matched nothing (2). The second group shrank from eight to two under the symbol-resolution changes described below; the recurring pattern in what remains is a parameter written with a subscript suffix, $c_i$, against a dataset field \texttt{cost}, where stripping the suffix before matching would succeed. The last four failures are one nonlinearity, one genuine parse failure, and two parameter whose index sets could not be determined from its arity.

\subsection{Infeasibility}

Eight easy and four hard problems compiled, ran, and reported infeasible. Two of them agree with an infeasible reference and are correct. The rest divide between formulations whose reference reports a different status and problems whose reference gives no verifiable answer, and the first group are over-constrained formulations. The most common source, on inspection, is a model emitting both an aggregate constraint and its per-element expansion, so that a resource limit is applied twice at different granularities. Analyzer pass 10 merges paired upper and lower bounds on the same expression but does not detect a constraint implied by a family of others, which would require a redundancy check we do not perform. This is an area of development and actively being worked on.

\subsection{Objective mismatches}

Seven problems solved whose objective disagreed with the reference. These are the most serious failures, because the system returns a confident answer that is wrong. The disagreements are gross, not numerical.

Median relative error is 0.54 on the easy split and 2.15 on the hard split, and no mismatch on either split falls within 0.1\% of the reference. Nothing here is a tolerance artifact; every mismatch is a structurally different model. Some have a reference objective of exactly zero, which is the signature of an objective built over the wrong term: the model optimizes a quantity that is identically zero at the optimum of the true problem. A further group, identified during earlier analysis, declares variables integer where the reference treats them as continuous, which tightens the feasible set and moves the optimum. Both are formulation errors that the compiler faithfully translates. Analyzer pass 7 rejects nonlinearity and pass 6 rejects undeclared symbols, but no pass can detect that a syntactically valid objective is the wrong objective.

This bounds what the deterministic approach can deliver. A compiler that faithfully translates a wrong formulation produces a wrong model, and problems in this category are wrong for that reason. Determinism makes the error reproducible and traceable; it does not make it go away.

\section{Determinism and regression testing}
\label{sec:determinism}

The emitter's output is a pure function of the IR, which makes byte-level regression testing possible. The compiler is covered by a six-gate matrix over a stored problem library, each gate a separate parametrized test: G1 compilation succeeds, G2 output parses as Python, G3 output executes, G4 the resulting model has the expected structure, G5 the solution matches a stored reference, and G6 the emitted source is byte-identical to a stored snapshot.  G6 is the gate that a generative pipeline cannot have. It asserts string equality against a checked-in file for the problem library, so any change in emitted code like a renamed variable, a reordered constraint, a different loop bound fails the test and has to be either fixed or explicitly re-baselined. During the compiler work described here that gate caught several changes whose effect on generated code was broader than intended. A pipeline that generates code by sampling cannot assert this, and has to fall back on testing that the model solves correctly, which is a far weaker condition and requires solving on every test run. The parser carries a further 49 unit tests and the end-to-end path 19 tests. This allows confidence in the output code reliability.

\section{Limitations}
\label{sec:limitations}

\paragraph{No detailed head-to-head baseline} The OptiMUS numbers we could verify from \cite{optimus2024} only discuss compilation with a commercial solver, and are not directly comparable to bulk of the numbers in Section~\ref{sec:results} other than accuracy results. Their solution accuracy is 80.3\% in the best case with o1, and 73.7\% with GPT-4o. The figure of ours that is closest in construction is 93.4\%, the match rate over all 320 problems with a numeric reference, rather than the 97.7\% of Table~\ref{tab:accuracy}, which conditions on the system having produced an answer. This limitation prevents a detailed comparison.

\paragraph{Formulation quality dominates the errors} Of the 25 compile failures, 15 are undeclared symbols. Of the mismatches, all inspected cases are structurally wrong models. Of the remaining 10 infeasibilities, the common cause is duplicated constraints. Most of the remaining error is upstream of the compiler, so compiler improvements have limited headroom against this benchmark; perhaps 4 of the 25 compile failures are addressable by better symbol resolution, on the evidence of Table~\ref{tab:compilefail}. The trajectory supports that reading: three changes to the parser and analyzer between our first measurement and this one moved the compile count from 325 to 329 without regressing anything, which is real but small.

\paragraph{Scope} The compiler targets LP and MILP. The seven case studies are excluded, quadratic objectives are rejected by an analyzer pass, and we make no claims about nonlinear, conic, or stochastic programs.

\paragraph{Benchmark provenance} NLP4LP problems come from textbooks and lecture notes, and we have not audited whether they appear in the training data of the model producing the formulations. This risk applies to the formulation stage; the compiler itself has no training data.

\section{Conclusion}

We show that separating the problem formulation stage from the translation stage costs little and buys a property that is otherwise unavailable. On 354 NLP4LP problems, a hand-written deterministic compiler produced solver code for 329 without consulting a language model, at a median compile cost below 10\,ms, with byte-identical output across runs. Of the 306 problems that could be evaluated, 299 matched the reference. The residual error is mostly outside the compiler. Fifteen of 25 compile failures are formulations referring to symbols they never declared, and all inspected objective mismatches are models that are structurally wrong rather than numerically imprecise. This is a useful thing to be able to say, and it is a consequence of the split: because the translation stage is deterministic and its failures are typed, the errors it does not cause can be attributed elsewhere. In an end-to-end generative pipeline, a wrong answer has no such decomposition.

The results allow us to make a strong claim; not all stages in solving an optimization problem need a language model, and some stages (esp. translation) are better off without language models. The benefits of deterministic behavior in these stages outweighs the benefits of introducing a language model.

{\small
\setlength{\itemsep}{0pt}

}


\begin{thebibliography}{99}
\setlength{\itemsep}{2pt}

\bibitem{optimus2023}
A.~AhmadiTeshnizi, W.~Gao, and M.~Udell.
\newblock OptiMUS: Optimization Modeling Using MIP Solvers and Large Language Models.
\newblock \emph{arXiv:2310.06116}, 2023.

\bibitem{optimus2024}
A.~AhmadiTeshnizi, W.~Gao, H.~Brunborg, S.~Talaei, and M.~Udell.
\newblock OptiMUS: Scalable Optimization Modeling with (MI)LP Solvers and Large
Language Models.
\newblock \emph{International Conference on Machine Learning (ICML)}, 2024.

\bibitem{solverllm2025}
D.~Li, X.~Zhao, L.~Yu, Y.~Liu, W.~Cheng, Z.~Chen, Z.~Chen, F.~Chen, C.~Zhao, and H.~Chen.
\newblock SolverLLM: Leveraging Test-Time Scaling for Optimization Problem via LLM-Guided Search.
\newblock \emph{Advances in Neural Information Processing Systems (NeurIPS)},
2025.

\bibitem{alphaopt2025}
M.~Kong, A.~Qu, X.~Guo, W.~Ouyang, C.~Jiang, H.~Zheng, Y.~Ma, D.~Zhuang, Y.~Tang, J.~Li, S.~Wang, H.~Koutsopoulos, H.~Wang, C.~Wu, and J.~Zhao.
\newblock AlphaOPT: Formulating Optimization Programs with Self-Improving LLM Experience Library.
\newblock \emph{arXiv:2510.18428}, 2025.

\bibitem{llamoco2024}
Z.~Ma, H.~Guo, J.~Chen, G.~Peng, Z.~Cao, Y.~Ma, and Y.-J.~Gong.
\newblock LLaMoCo: Instruction Tuning of Large Language Models for Optimization Code Generation.
\newblock \emph{IEEE Transactions on Evolutionary Computation}, 2026.

\bibitem{li2023milp}
Q.~Li, L.~Zhang, and V.~Mak-Hau.
\newblock Synthesizing Mixed-Integer Linear Programming Models from Natural Language Descriptions.
\newblock \emph{arXiv:2311.15271}, 2023.

\bibitem{smilo2026}
Q.~Li, L.~Zhang, and V.~Mak-Hau.
\newblock An LLM-powered MILP modelling engine for workforce scheduling guided
by expert knowledge.
\newblock \emph{International Journal of Production Research}, 2026.
\newblock DOI: 10.1080/00207543.2026.2663386.

\bibitem{mamo2025}
X.~Huang, Q.~Shen, Y.~Hu, A.~Gao, and B.~Wang.
\newblock LLMs for Mathematical Modeling: Towards Bridging the Gap between
Natural and Mathematical Languages.
\newblock \emph{Findings of the ACL: NAACL 2025}, pp.~2678--2710, 2025.

\bibitem{component2025}
D.~Refai and M.~Ahmed.
\newblock Peering Inside the Black Box: Uncovering LLM Errors in Optimization Modelling through Component-Level Evaluation.
\newblock \emph{arXiv:2510.16943}, 2025.

\bibitem{optichat2025}
H.~Chen, G.~E.~Constante-Flores, K.~S.~I.~Mantri, S.~M.~Kompalli, A.~S.~Ahluwalia, and C.~Li.
\newblock OptiChat: Bridging Optimization Models and Practitioners with Large Language Models.
\newblock \emph{INFORMS Journal on Data Science.}, 2025.

\bibitem{survey_llmopt2024}
S.~Huang, K.~Yang, S.~Qi, and R.~Wang.
\newblock When Large Language Model Meets Optimization.
\newblock \emph{Swarm and Evolutionary Computation, 90, 101663.}, 2024.

\bibitem{survey_or2025}
Y.~Wang and K.~Li.
\newblock Large Language Models in Operations Research: Methods, Applications,
and Challenges.
\newblock \emph{arXiv:2509.18180}, 2025.

\bibitem{docp2024}
S.~Wasserkrug, L.~Boussioux, D.~den~Hertog, F.~Mirzazadeh, I.~Birbil, J.~Kurtz, and D.~Maragno.
\newblock From Large Language Models and Optimization to Decision Optimization CoPilot: A Research Manifesto.
\newblock \emph{arXiv:2402.16269}, 2024.

\bibitem{lm4opt2024}
T.~Ahmed and S.~Choudhury.
\newblock LM4OPT: Unveiling the Potential of Large Language Models in Formulating Mathematical Optimization Problems.
\newblock \emph{ INFOR: Information Systems and Operational Research, 62(4), 559-572.}, 2024.

\bibitem{nl4opt2023}
R.~Ramamonjison, T.~T.~Yu, R.~Li, et~al.
\newblock NL4Opt Competition: Formulating Optimization Problems Based on Their
Natural Language Descriptions.
\newblock \emph{NeurIPS Competition Track}, 2023.

\end{thebibliography}
\end{document}